\documentclass[twoside,11pt]{article}
\usepackage{amsmath} 
\usepackage[e]{esvect}
\usepackage{graphicx}
\usepackage[adobe-utopia]{mathdesign}
\usepackage[T1]{fontenc}
\usepackage[english]{babel}
\usepackage{diagbox,booktabs}
\usepackage{color}

\usepackage{hyperref}

\title{\textsc{A General Inequality for Packings of Boxes}}
\date{}
\author{{Krzysztof Przes{\l}awski}\\
}
\newtheorem{pr}{\sc Proposition}
\newtheorem{lemat}[pr]{\sc Lemma}
\newtheorem{tw}[pr]{\sc Theorem}

\newtheorem{df}{\sc Definition}

\newtheorem{uw}{\sc Remark}
\newtheorem{uwi}[uw]{\sc Remarks}
\newenvironment{uwa}{\begin{uw} \rm}{\end{uw}}

\newtheorem{nap}{\sc Example}

\newtheorem{nps}[nap]{\sc Examples}
\newenvironment{nprs}{\begin{nps}\rm}{\end{nps}}

\def\kal #1 #2{\mathscr{#1}^{#2}}
\def\proof{\noindent \textit{Proof.\,\,\,}}

\def\er{\mathbb{R}}

\def\Aut #1 #2{\operatorname{Aut}^{#1} (#2)}

\def\bred #1 {\colorbox{red}{ #1}}
\def\red #1 {{\color{red} #1 }}

\graphicspath{{figpdfkolor/}}
\DeclareGraphicsExtensions{.pdf,.png} 

\begin{document}
\maketitle
\begin{abstract}
Keller packings and tilings  of  boxes are investigated.  Certain general inequality measuring a complexity of such systems is proved.  A straightforward application 
to the unit cube tilings is given:

Let $\boldsymbol m=(m_1, \ldots, m_d)$ be a vector with integer coordinates such that $m_i\ge 2$ for every $i$. Let $\boldsymbol m\mathbb Z^d=m_1\mathbb Z\times\cdots\times m_d\mathbb Z$.  We identify the $d$-dimensional torus $\mathbb R^d/ {{\boldsymbol m}\mathbb Z^d}$ with  $\mathbb T_{\boldsymbol m}^d=[0,m_1)\times[0,m_2)\cdots\times [0,m_d)$. The addition $\oplus$ in $\mathbb T_{\boldsymbol m}^d$ is  transferred from $\mathbb R^d/\boldsymbol m\mathbb Z^d$; that is, $x\oplus y\in \mathbb T_{\boldsymbol m}^d$ is the only element satisfying the congruence $x\oplus y \equiv x+y\pmod{\boldsymbol m}$. A family $[0,1)^d\oplus T=\{[0,1)^d\oplus t\colon t \in T\}$ of translates of the unit cube $[0,1)^d$ is said to be a \textit{cube tiling} of $\mathbb T^d_{\boldsymbol m}$ if it is a partition of $\mathbb T^d_{\boldsymbol m}$.
Let  us set
$$
p_i(T)=\{ t_i\bmod 1\colon t \in T\}  
$$
for $i\in \{1,\ldots, d\}$. 
A  straightforward consequence of our inequality reads that if $\boldsymbol m=(n,n,\ldots,n)$, then 
$$
\sum_i|p_i(T)|\le \frac {n^d-1}{n-1}
$$ 
The inequality is tight. The extremal cube tilings for which it becomes an equality are classified.  
\bigskip

\end{abstract}

\section{Introduction} 

Let $\boldsymbol m=(m_1, \ldots, m_d)$ be a vector with integer coordinates such that $m_i\ge 2$ for every $i$. Let $\boldsymbol m\mathbb Z^d=m_1\mathbb Z\times\cdots\times m_d\mathbb Z$. We identify the $d$-dimensional torus $\mathbb R^d/ {{\boldsymbol m}\mathbb Z^d}$ with  $\mathbb T_{\boldsymbol m}^d=[0,m_1)\times[0,m_2)\cdots\times [0,m_d)$. The addition $\oplus$ in $\mathbb T_{\boldsymbol m}^d$ is  transferred from $\mathbb R^d/\boldsymbol m\mathbb Z^d$; that is, $x\oplus y\in \mathbb T_{\boldsymbol m}^d$ is the only element satisfying the congruence $x\oplus y \equiv x+y\pmod{\boldsymbol m}$. A family $[0,1)^d\oplus T=\{[0,1)^d\oplus t\colon t \in T\}$ of translates of the unit cube $[0,1)^d$ is said to be a \textit{cube tiling} of $\mathbb T^d_{\boldsymbol m}$ if it is a partition of $\mathbb T^d_{\boldsymbol m}$.
Let  us define the parameters of  this tiling:
$$
p_i(T)=\{ t_i\bmod 1\colon t \in T\}  
$$
for $i\in \{1,\ldots, d\}$. The following problem can be derived  from Dutour and Itoh \cite[Conjecture 5.4]{DI}:
\begin{quotation}
Let  ${\boldsymbol m}=(2,2,\ldots, 2)$.  Show that    $\sum_i|p_i(T)|\le 2^d-1$. Moreover, prove that the equality takes place only in the case of tilings that are obtained by a special lamination construction.
\end{quotation}

One of our goals  is to show that for any integer $n\ge 2$, if $\boldsymbol m=(n,n,\ldots,n)$, then 
$$
\sum_i|p_i(T)|\le \frac {n^d-1}{n-1}
$$ 
and the inequality is tight. The extremal cube tilings for which this inequality becomes an equality are classified (Theorem \ref{C}).   Theorem \ref{C} appears to be a consequence of a more general result concerning the so-called Keller packings and partitions (Theorem \ref{B}).       

\section{Keller families}
Cube tilings have attracted a great deal of attention lasting for over one hundred years mainly in connection with Minkowski's and  Keller's  conjectures, e. g., \cite{CSz, D, H, Ke1, Ke2, KL, KP2, Ko, Ma, M, LP1, LP2, LS1, LS2, P, SS, Sz, Z}.  They are also related to Fuglede's conjecture \cite{F}. In particular, it appeared quite recently that they  are in  a 1--1 correspondence  with the so-called exponential bases  in $L^2([0,1]^d)$, e.g., \cite{IP, LRW}. 
  
Kearnes and Kiss \cite{KK} stated the following combinatorial problem, which, as they mentioned, was completely unrelated to algebraic questions they were considering:
\begin{quotation}
``\textbf{Problem 5.5.} Let $A_1,\ldots, A_k$ be nonempty sets. If the rectangle $A = A_1\times\cdots \times A_k$
is partitioned into less than $2^k$ rectangular subsets, does it follow that there exists a
rectangular subset in this partition which has full extent in some direction $i$?''
\end{quotation}
This problem is settled in \cite{ABHK}.   It is proved that a minimal partition of the rectangle $A$ into rectangular subsets without full extent in any direction has to have $2^k$ elements. Such partitions are thoroughly investigated in \cite{GKP, KP1} .  It appears that they are combinatorially equivalent to unit cube tilings of the torus $\mathbb T^k_{\boldsymbol m}$, where $\boldsymbol m=(2,\ldots, 2)$.  This observation suggest that even if one investigates unit cube tilings it can be profitable  to work with more abstract structures suggested by Problem 5.5.  We follow this suggestion in the present paper.   

Let $[0,1)^d + T=\{[0,1)^d+t\colon t \in T\}$ be a cube tiling of $\mathbb R^d$. As observed by Keller \cite{Ke1}, for every pair $s, t\in T$ of different vectors there is $i$ such that $s_i-t_i$ is a non-zero integer (see also \cite{IP}). Obviously, there is a counterpart of this result for cube tilings of tori, since they can be interpreted as periodic tilings of $\mathbb R^d$. A thorough inspection shows that Keller's result is a consequence of the fact that if we take two unit segment  partitions $[0,1)+S$ and $[0,1)+T$ of $\mathbb R$, then $T=S$ or for every pair of proper subsets $P\subset S$, $Q\subset T$ the union of $[0,1)+P$ does not coincide with the union of $[0,1)+Q$. This  observation motivates our further discussion. 
\subsection{Partitions}
A partition $\pi$ of a set $U$ is \textit{non-trivial}, if $\pi \neq \{U\}$. Two partititions $\pi$ and $\pi'$ are \textit{independent} if $\{U\}=\pi\vee\pi'$; that  is, $\{U\}$ is  a unique partition coarser than these two. Consequently,

\begin{lemat} 
\label{C forte}
If partitions $\pi$ and $\pi'$ of $U$  are independent, $\rho$ is a proper subfamily of $\pi$ and $\rho'$ is a proper subfamily of $\pi'$, then the union of $\rho$ and that of $\rho'$ do not coincide.
\end{lemat}

In particular, it follows from Lemma \ref{C forte} that if $\Pi$ is a family of pairwise independent partitions of $U$, then for any set $A$ belonging to the union of $\Pi$ there is a unique partition $\pi$ such that $A\in \pi$. 
  
A family $\Pi$ of  partitions of $U$ is \textit{complete} if every partition $\rho$ of $U$ whose parts belong to the union of $\Pi$ is a memeber of $\Pi$.  It is \textit{unital} if $\{U\}$ is its member.  
\begin{nprs}
\label{np}
\medskip
($\alpha_1$)  Let $U$ be a non-empty set. Then $\Pi=\{\{A, U\setminus A\} \colon A\subseteq U,\, A\not\in\{U, \emptyset\}\}\cup \{\{U\}\}$ is pairwise independent and unital.  In general, it 
is not complete.  

\medskip
\noindent ($\alpha_2$) This example has been already mentioned. Let $U=\mathbb T_m=[0,m)$, where $m$ is a positive integer. For every $t\in [0,1)$, let us define $\pi_t= \{[0,1)\oplus t\oplus k\colon k \in \{0,\ldots m-1\}\}$. Then $\Pi=\{\pi_t\colon t \in [0,1)\}\cup \{\{U\}\}$ is a complete and unital family of pairwise independent partitions of $U$. The family of all non-trivial partitions belonging to $\Pi$ is simply the family of all unit segment tilings of $U$.

\medskip
\noindent ($\alpha_3$) Let us fix a hexagon $H$ in $\mathbb R^2$ with three consecutive edges removed. Clearly $\mathbb R^2$ can be tiled by translates of $H$. In fact, each such a tiling is a partition of $\mathbb R^2$. Let us assume that $\pi$ is a member of $\Pi$ if and only if  $\pi$ is a tiling of $\mathbb R^2$ by translates of  $H$ or $\pi$ is equal to a trivial partition of $\mathbb R^2$. Again $\Pi$ is complete, unital and consists of pairwise inedependent paritions. Tilings discussed here are infinite.  Since they are periodical, one can define corresponding finite tilings of two-dimensional tori by translates of $H$.            
\end{nprs}

\subsection{Boxes. Keller's condition}
Let $X_i$, $i\in\{1, \ldots,d\}$, be a sequence of sets such that each of them has at least two elements. Let $\Pi_i$ be a finite family of pairwise independent finite partitions of $X_i$. Let $X=X_1\times\cdots\times X_d$ and $\mathscr U_i$ be the union of $\Pi_i$. Let us define the family of boxes  $\mathscr B=\mathscr U_1\otimes\cdots\otimes\mathscr U_d$: 
$$ 
K=K_1\times\cdots\times K_d  \in \mathscr B\quad  \text{if  and only if}\quad   K_1\in \mathscr U_1, \ldots, K_d\in \mathscr U_d.
$$
We say that a box $K\in \mathscr B$ is \textit{proper}, if  $K_i\neq X_i$, for every $i$. Otherwise, $K$ is \textit{improper}.   The triple $(X, \Pi, \mathscr B)$ is called a \textit{system of $d$-boxes}.  Occasionally,  we abuse slightly the terminology referring to $\mathscr B$ alone as the system of boxes.   

Now, we demonstrate an abstract version of Keller's result concerning unit-cube tilings of tori. 
\begin{pr}
\label{kellerian}
Let $(X, \Pi, \mathscr B)$ be a system of $d$-boxes. Let $\mathscr G\subseteq \mathscr B$ be a partition of $X$.  If   $\Pi_i$ is complete for every $i\in\{1,\ldots, d\}$, then $\mathscr G$ satisfies {\em Keller's condition}:
\begin{quote}
If    $K$, $L\in \mathscr G$ and $K\neq L$, then there is  $i$ such that  $K_i$, $L_i$ are different parts of a partition $\pi\in \Pi_i$.
\end{quote} 
\end{pr}

\proof Clearly, we may assume that each $\Pi_i$ is unital.  Suppose now that our proposition is false. Then we can choose such a partition $\mathscr G$, among those violating Keller's condition, that there is a pair $K, L\in \mathscr G$ so that $\{K,L\} $ violates Keller's condition and   the set $I=\{i\colon K_i\cap L_i=\emptyset\}$ is minimal in the sense of inclusion. Let us pick any $s\in I$. Let $\pi\in \Pi_s$ be chosen so that 
$L_s\in \pi$. Let as define a `pile' $\mathscr P=\{G\in \mathscr G\colon G_s\in \pi\}$. Let us fix any $x=(x_1,\ldots, x_d)$ belonging to a member $H$ of $\mathscr P$. Let 
$$
\ell = \{x_1\}\times\cdots\times\{x_{s-1}\}\times X_s \times \{x_{s+1}\}\times\cdots \times \{x_d\}.
$$
Let $\mathscr G_\ell=\{G\in \mathscr G\colon G\cap \ell \neq\emptyset\}$. Since the sets $G\cap \ell$, for $G\in \mathscr G_\ell$ form a partiton of $\ell$, the sets  $(G\cap \ell)_s = G_s$, for $G \in \mathscr G_\ell$, form a partiton $\rho$ of $X_s$. As $H\in \mathscr G_\ell$, it follows that $\pi$ and $\rho$ share the same part $H_s$. Since $\Pi_s$ is complete, $\rho$ has to belong belong to $\Pi_s$. As $\pi$ and $\rho$ are not independent, they have to coincide. Let $P$ be the union of $\mathscr P$. It is clear that $P=\bigcup_{x\in P} \ell_x$. Therefore, we can imagine $P$ as a kind of cylinder. For every box $B\in \mathscr B$ we can define a new box $B^s$ by replacing factor $B_s$ by $X_s$. Let $\mathscr H=\{ G^s\colon \text{$G\in \mathscr P$ and $G_s=L_s$}\}$. The union of $\mathscr H$ coincides with $P$. Therefore, $(\mathscr G\setminus \mathscr P) \cup \mathscr H$ is a partition of $X$. Moreover $L^s$ and $K$ are different and $J=\{i \colon K_i\cap L^s_i=\emptyset\}=I\setminus\{s\}$, which contradicts the minimality of $I$. $\square$       

Let $\mathscr G$ be a non-empty subfamily of  a system of  boxes $\mathscr B$. We call it \textit{Keller family} if it satisfies Keller's condition. Clearly, Keller families consist of disjoint boxes, therefore, one can find justified to refer to them as \textit{Keller packings}.  

A set $G\subseteq X $ which is the union of a  Keller family $\mathscr G\subseteq \mathscr B$ is called a \textit{polybox}
while $\mathscr G$ itself is called a \textit{suit} for $G$. A polybox can have more than one suit.  
The suit $\mathscr G$ is \textit{proper} if each box $K\in \mathscr G$ is proper; otherwise, it is \textit{improper}.  

\begin{uwa}
Systems of $d$-boxes $(X, \Pi,\mathscr B)$ such that $(\Pi_i\colon i=1,\ldots,d)$ are as described in Example \ref{np}.$\alpha_1$ are investigated in \cite{GKP, KP1}.
\end{uwa}

\section{Main results}

Let $(X,\Pi, \mathscr B)$ be a system of $d$-boxes. Let $i \in \{1,2\ldots, d\}$ and let $\mathscr V$ be a nonempty subfamily of  $\mathscr U_i=\bigcup \Pi_i$.  For every $\mathscr G\subseteq \mathscr B$, we define the \textit{restriction} $\mathscr G|\mathscr V$ of $\mathscr G$ with respect to $\mathscr V$ by the equality
$$
\mathscr G|\mathscr V=\{K\in \mathscr G\colon K_i\in \mathscr V\}.
$$
If $\mathscr V$ is a singleton of $A$, then  we write $\mathscr G|A$ rather than $\mathscr G|\{A\}$.   

Let $\pi \in \Pi_i$ be a non-trivial partition. A Keller family $\mathscr G$ is \textit{laminated with respect to $\pi$ if $\mathscr G|\pi$}  coincides with $\mathscr G$. It is \textit{$\{i\}$-laminated} if it is laminated with respect to a non-trivial partition from $\Pi_i$. Finally, the family $\mathscr G$ is \textit{laminated} if  it is $j$-laminated for some $j\in\{1,\ldots, d\}$. 
   

A polybox $W\subseteq X$ is said to be  an $i$-\textit{cylinder} or  shortly  a \textit{cylinder}   if for exery $x\in W$ the \textit{line}
$
\ell= \{x_1\}\times\cdots\times\{x_{i-1}\}\times X_i \times \{x_{i+1}\}\times\cdots \times \{x_d\}
$
is contained in $W$. 

A Keller family $\mathscr C$ is said to be an $i$-\textit{pile} or  shortly \textit{pile}  if it is a suit for an $i$-cylinder $W$ and is laminated with respect to a non-trivial partition $\pi \in \Pi_i$. For every  $A\in \pi$, let us define

$$
\mathscr C^A=\{K_1\times\cdots \times K_{i-1}\times X_i\times K_{i+1}\times\cdots\times K_d\colon K\in \mathscr C|A\}
$$

We call $\mathscr C^A$ an \textit{elementary aggregate} of the pile $\mathscr C$ (with respect to $A$). It is clear that $\mathscr C^A$ is a suit  for $W$. 

Let $i \in\{1,\ldots,d\}$. Let $\pi\in \Pi_i$ be a non-trivial partition,  and $\mathscr G\subseteq \mathscr B$ be a Keller family.  We write that $\pi$ is $\textit{present}$ in $\mathscr G$ if there is $K\in \mathscr G$ such that $K_i\in \pi$.  
Partition $\pi$ is $\textit{hidden}$ in $\mathscr G$ if $\mathscr G|\pi$ is  a suit  for an  $i$-cylinder.  
Moreover, $\pi$ is said to be $\textit{exposed}$ in $\mathscr G$ if it is present and  not hidden in $\mathscr G$.  

We denote by $C_i(\mathscr G)$ the family of all partitions $\pi\in \Pi_i$ which are hidden in $\mathscr G$.  Let 
$$
c_i(\mathscr G)= \sum_{\pi \in C_i(\mathscr G)} (|\pi|-1)\quad \text{and}\quad  c(\mathscr G)=\sum_{i=1}^d c_i(\mathscr G).
$$

We shall need a recursive definition of a \textit{multipile}. Let $\mathscr G\subseteq\mathscr B$ be a Keller family. Suppose that one of the following two conditions is satisfied:
\begin{enumerate}
\item there is a box $K\in \mathscr B$ such that  $\mathscr G = \{K\}$; 
\item there is $i$ and non-trivial $\pi \in \Pi_i$ such that: 
\begin{itemize}
\item[-] $\mathscr G$ is an $i$-pile laminated with respect to $\pi$,
\item[-] for every $A\in \pi$, $\mathscr G|A$ is a multipile,
\item[-] for every $k\neq i$ and every $A, B\in \pi$, if $A\neq B$, then $C_k(G|A)\cap C_k(G|B)=\emptyset$.
\end{itemize}
\end{enumerate}
Then $\mathscr G$ is a multipile. 
\begin{uwa}
\label{M}
It follows easily by induction that  each multipile is a suit for a box. 
\end{uwa}

\begin{tw}
\label{B} 
Let $(X, \Pi, \mathscr B)$ be a system of boxes. Suppose $\mathscr G\subseteq  \mathscr B$ is a Keller family. Then 
$$
c(\mathscr G)\le |\mathscr G|-1.
$$
The inequality becomes an equality if and only if  $\mathscr G$ is a multipile.
\end{tw} 

We precede the proof of the theorem by two lemmas. 

\begin{lemat} 
\label{A}
Let $(X, \Pi, \mathscr B)$ be a system of boxes. Let $\mathscr G, \mathscr G'\subseteq \mathscr B$ be two suits for the same  polybox $G$. If  $\pi\in \Pi_i$  is exposed in $\mathscr G$, then it is exposed in $\mathscr G'$.
\end{lemat}

\proof Since $\pi$ is exposed, $\mathscr G|\pi$ is non-empty, and  is not an $i$-pile. Therefore, there is  $u \in X$ such that the line
$
\ell=\{u_1\}\times\cdots\{u_{i-1}\} \times X_i\times \{u_{i+1}\}\times \cdots\times\{u_d\} 
$
intersects some of the members of the family $\mathscr G|\pi$ but is not covered by this family. Let $\rho=\{A_i\colon \text{$A \in \mathscr G$ and $A\cap \ell \neq \emptyset$}\}$. Since $\mathscr G$ is a Keller family, all elements of $\rho$ have to be parts of one partition belonging $\Pi_i$. Now, since $\pi$ and $\rho$ have elements in common and $\Pi$ is pairwise independent,  $\rho\subset \pi$. Let $R$ be the  union of $\rho$. Then         
$
G\cap \ell= \{u_1\}\times\cdots\{u_{i-1}\} \times R\times \{u_{i+1}\}\times \cdots\times\{u_d\}.
$
Let $\rho'=\{A_i\colon \text{$A \in \mathscr G'$ and $A\cap \ell \neq \emptyset$}\}$. Since $\mathscr G'$ is a suit for $G$, we deduce that the union of $\rho'$ is $R$. By Lemma \ref{C forte}, it follows that $\rho'\subseteq \pi$. Therefore $\pi$ is exposed in $\mathscr G'$.     

As an immediate corollary we have:
\begin{lemat}
\label{H}
Let $(X, \Pi, \mathscr B)$ be a system of boxes. Let $\mathscr G, \mathscr G'\subseteq \mathscr B$ be two suits for the same  polybox $G$. If $\pi \in \Pi_i$ is hidden in $\mathscr G$ and present in $\mathscr G'$, then it is hidden in $\mathscr G'$.  
\end{lemat}

\noindent\textit{Proof of the theorem.}
We shall proceed by induction with respect to the cardinality of $\mathscr G$.

Suppose that $\pi\in \Pi_i$ is hidden in $\mathscr G$. Then $\mathscr C=\mathscr G|\pi$ is an $i$-pile. Fix an $A\in \pi$. Then $\mathscr C$ and $\mathscr C^A$  are suits for the same $i$-cylinder.  Observe that  the family $\mathscr G'= (\mathscr G\setminus \mathscr C)\cup \mathscr C^A$ is kellerian.  To this end it suffices to pick $K\in  \mathscr G\setminus \mathscr C$ and $L \in \mathscr C^A$, and show that $\{K,L\}$ is kellerian.  By the definition of $\mathscr C^A$, the box $L'=L_1\times \cdots \times L_{i-1}\times A_i\times L_{i+1}\times \cdots\times L_d$ is a member of   $\mathscr C$. Therefore, $\{K, L'\}$ is kellerian as a subfamily of $\mathscr G$.    It follows from the definition of $\mathscr C$that  $K_i\not \in \pi$ while $L'_i=A\in \pi$.  Consequently, there is $j \neq i$ such that $K_j$ and $L'_j=L_j$ are different parts of the same partition $\tau\in \Pi_j$, which proves that $\{K,L\}$ is kellerian as expected.  It is clear now that $\mathscr G$ and $\mathscr G'$ are suits for the same polybox. 

As an immediate consequence of the definition of $\mathscr G'$ one has

\begin{equation}
\label{P1}
C_i(\mathscr G)=C_i(\mathscr G')\,\dot{\cup}\, \{\pi\}, 
\end{equation}

where the $\dot{\cup}$ symbol denotes the disjoint union.  

Suppose now that $k\neq i$ and $\rho\in \Pi_k$ is hidden in $\mathscr G$. If $\rho$ is present in $\mathscr G'$, then, by Lemma \ref{H}, it is hidden in $\mathscr G'$. If it is not present in $\mathscr G'$, then there is such a $B\in \pi\setminus\{A\}$  that $\rho$ is present in $\mathscr G|B$, which is equivalent to saying that $\rho$ is present in $\mathscr C^B$. Let $\mathscr G''= (\mathscr G\setminus \mathscr C)\cup \mathscr C^B$. It is clear now that $\rho$ is hidden in $\mathscr G''$. Since $\rho$ is not present in $\mathscr G$, it cannot be present in $\mathscr G\setminus \mathscr C$. Consequently, $\rho$ is hidden $\mathscr C^B$, which readily implies that it is hidden in  $\mathscr G|B$. Therefore, we have just shown that for every $k\neq i$

\begin{equation}
\label{P2}
C_k(\mathscr G)= C_k(\mathscr G')\cup \bigcup_{B\in \pi\setminus\{A\}} C_k(\mathscr G|B). 
\end{equation}

By (\ref{P1}) we have,
\begin{equation}
\label{c1}
c_i(\mathscr G)=c_i(\mathscr G')+|\pi|-1. 
\end{equation}
By (\ref{P2}),
\begin{equation}
\label{c2}
c_k(\mathscr G)\le c_k(\mathscr G')+\sum_{B\in \pi\setminus\{A\}} c_k(\mathscr G|B). 
\end{equation}
Consequently,
\begin{equation}
\label{c3}
c(\mathscr G)\le c_i(\mathscr G')+|\pi|-1+ \sum_{k\neq i}\left(c_k(\mathscr G')+\sum_{B\in \pi\setminus\{A\}} c_k(\mathscr G|B)\right). 
\end{equation}
Since $c_i(\mathscr G|B)=0$, we get
\begin{equation}
\label{c4}
c(\mathscr G)\le c(\mathscr G') +|\pi|-1+\sum_{B\in \pi\setminus\{A\}} c(\mathscr G|B)). 
\end{equation}
As $|\mathscr G'|< |\mathscr G|$ and $|\mathscr G|B|< |\mathscr G|$, we conclude by induction 
\begin{equation}
\label{c5}
c(\mathscr G)\le |\mathscr G'|-1 +|\pi|-1 +\sum_{B\in \pi\setminus\{A\}} (|\mathscr G|B|-1)=|\mathscr G'|-1 +\sum_{B\in \pi\setminus\{A\}} |\mathscr G|B|.
\end{equation}
By the definition of $\mathscr G'$, we have $|\mathscr G'|=|\mathscr G\setminus \mathscr C|+|\mathscr C^A|$. Since $|\mathscr C^A|=|\mathscr G|A|$, we obtain
\begin{equation}
\label{c6}
c(\mathscr G)\le |\mathscr G\setminus \mathscr C|-1+\sum_{B\in \pi}|\mathscr G|B|=|\mathscr G\setminus \mathscr C|-1+|\mathscr C|=|\mathscr G|-1. 
\end{equation}

It remains to discuss the case of equality. The fact that the equality holds for each multipile is an easy part of the proof, and as such  is left to the reader.

If $c(\mathscr G)=|\mathscr G|-1$, then all the above inequalities become equalities. In particular,  (\ref{c2}) becomes an equality. This equality combined with (\ref{P2}) readily implies that the sets $C_k(\mathscr G|B)$, for $B\in \pi \setminus\{A\}$ are pairwise disjoint. Since $A$ is arbitrary, we conclude that 

\begin{itemize}
\item[($\beta$)] For every $k\neq i$, the indexed family $(C_k(\mathscr G|B)\colon B\in \pi)$ consists of pairwise disjoint members. 
\end{itemize}

As (\ref{c4}) and (\ref{c5}) become equalities and $A$ can be replaced by any element of $\pi$, it follows that $c(\mathscr G|B)=|\mathscr G|B|-1$ for $B\in \pi$. By induction hypothesis, the sets $\mathscr G|B$ are multipiles. This fact combined with ($\beta$) and the definition of $\mathscr C$ leads to the conclusion that $\mathscr C=\mathscr G|\pi$ is a multipile. Thus, in the case $\mathscr G=\mathscr C$ there is nothing more to prove. Suppose now that $\mathscr G$ is different from $\mathscr C$. By Remark \ref{M}, $\mathscr C$ is a suit for a box $K\in \mathscr B$. Let us set $\mathscr G'''= (\mathscr G\setminus \mathscr C)\cup \{K\}$.  Clearly, $\mathscr G'''$ and $\mathscr G$ are suits for the same polybox. Suppose that $ \rho \in C_k(\mathscr G)$ for a certain $k\neq i$. If $\rho$ is present in $\mathscr G\setminus \mathscr C$, then by Lemma \ref{H} it is hidden in $\mathscr G'''$. Clearly, $C_k(\mathscr G''')\subseteq C_k(\mathscr G)$. If $\rho$ is not present in $C_k(\mathscr G\setminus \mathscr C)$, then it has to be hidden in some of  the sets $\mathscr G|B$, where $B\in \pi$. Consequently, 
\begin{equation}
\label{P3}
C_k(\mathscr G)= C_k(\mathscr G''')\cup \bigcup_{B\in \pi} \mathscr G|B. 
\end{equation}
Suppose now that $ \rho \in C_k(\mathscr G''')$. Since $\pi$ is nontrivial, $|\pi|>1$. By ($\beta$), there is $A\in \pi$ such that $\rho\not \in C_k(\mathscr G|A)$. Since (\ref{c2}) becomes an equality for this paricular $A$ as for any other, the right-hand side expression of the equality (\ref{P2}) consists of disjoint sets. Therefore, $\rho$ cannot belong to any of the sets $C_k(\mathscr G|B)$, $B\in \pi$. Thus, the right-hand side expression of the equality (\ref{P3}) consists of disjoint sets. As an immediate consequence,
$$
\sum_{k\neq i} c_k(\mathscr G)= \left(\sum_{k\neq i} c_k(\mathscr G''')\right)+ \sum_{B\in \pi} c(\mathscr G|B).
$$
As we already know, $c(\mathscr G)=|\mathscr G|-1$ implies $c(\mathscr G|B)=|\mathscr G|B|-1$. Moreover, it follows from the definition of $\mathscr G'''$ that  $c_i(\mathscr G)=c_i(\mathscr G''')+(|\pi|-1)$. Thus,
$$
|\mathscr G|-1= c(\mathscr G''')+|\pi|-1 +\sum_{B\in \pi} (|\mathscr G|B|-1)= c(\mathscr G''') + |\mathscr C|-1.
$$
And
$$
c(\mathscr G''')= |\mathscr G\setminus \mathscr C|=|\mathscr G'''|-1.
$$
Since $|\mathscr G'''| < |\mathscr G|$, we deduce by induction that $\mathscr G'''$ is a multipile. In particular, it is a $j$-pile for certain $j\neq i$. Then $\mathscr G$ is a $j$-pile. As we have already discussed the case where $\mathscr G$ is a pile, the proof is complete.    $\Box$  

\begin{uwa}
\label{podzial}
If $\mathscr G\subseteq B$ is a partition of $X$, then, by  Proposition \ref{kellerian}, it is automatically Kellerian.  In particular, each unit cube tiling $[0,1)^d+T$ of $\mathbb T^d_{\boldsymbol m}$ is Kellerian.  
\end{uwa}
Now we are ready to prove our main result concerning cube tilings.

\begin{tw} 
\label{C}
Let $[0,1)^d\oplus T=\{[0,1)^d\oplus t\colon t \in T\}$ be a cube tiling of $\mathbb T^d_{\boldsymbol m}$. Let 
\begin{equation}
\label{P}
p_i(T)=\{ t_i\bmod 1\colon t \in T\}  
\end{equation}
for $i\in \{1,\ldots, d\}$.
If ${\boldsymbol m}=(n,n,\ldots, n)$ for an integer $n\ge 2$, then 
\begin{equation}
\label{eb}
p(T):=\sum_i|p_i(T)|\le \frac {n^d-1}{n-1}. 
\end{equation}
The equality takes place if and only if $[0,1)^d\oplus T$ is a multipile.
\end{tw}
\proof
Let $X_i=\mathbb T_m = [0,m)$ for  $i \in \{1,\ldots, d\}$ . Let  $\Pi_i$  be as described in Example \ref{np}.$\alpha_2$.  Let $(X, 
\Pi, \mathscr B)$ be the corresponding system of $d$-boxes. Let $\mathscr G=[0,1)^d\oplus T$. Clearly, $\mathscr G\subset \mathscr B$.  By Proposition 2, $\mathscr G$ satisfies Keller's condition. Therefore,  Theorem \ref{B}  applies to $\mathscr G$.    

For every $i \in \{1,\ldots, d\}$, we have $c_i(\mathscr G)=(n-1)|p_i(T)|$. Moreover, $|\mathscr G|=n^d$. The inequality (\ref{eb}) follows as an immediate consequence of Theorem \ref{B}, so does the case of equality. $\square$     

\medskip
An abstract version of Theorem \ref{C} holds true as well.  

\begin{tw}
\label{D}
Let $(X, \Pi, \mathscr B)$ be a system of $d$-boxes. Suppose there is  an integer $n$ such that for every $i\in\{1,\ldots, d\}$ all nontrivial partitions belonging to $\Pi_i$ are of cardinality $n$. Suppose $\mathscr G\subset \mathscr B$ is a Keller partition of $X$ which consists of proper boxes. Then 
$$
\sum_i |C_i(\mathscr G|\le \frac{n^d-1}{n-1}. 
$$
The equality takes place if and only if  $\mathscr G$ is a multipile.
\end{tw} 

Clearly, the proof reduces to showing that 
\begin{equation}
\label{N}
|\mathscr G|=n^d.
\end{equation} 
We begin by discussing a certain general construction, which is instrumental in the theory of Keller packings.

We may assume  that each $\Pi_i$,  $i\in \{1,\ldots, d\}$, is unital.
Let $Y_i=\prod_{\pi\in \Pi_i}\pi$, for $i \in \{1,\ldots, d\}$.  Let $A\in \mathscr U_i$, where as before $\mathscr U_i=\bigcup \Pi_i$. Since $\Pi_i$ is a family of pairwise independent partitions, there is exactly one $\rho \in \Pi_i$ such that $A\in \rho$. Let us define $\hat{A}=\prod_{\pi\in \Pi_i} \hat{A}_\pi \subseteq Y_i$ by the formula
$$
\hat{A}_\pi=
\left\{ \begin{array}{rl}
\pi, &\text{if $\pi \neq \rho$,}\\
\{A\}, & \text{if $\pi=\rho$.}
\end{array}
\right.
$$
Observe that for every partition $\pi\in \Pi_i$, the family $\hat{\pi}=\{\hat{A}\colon \text{for $A\in\pi$}\}$ is a partition of $Y_i$. Moreover, the family of  partitions $\hat{\Pi}_i=\{\hat{\pi}\colon \pi \in \Pi_i\}$ is pairwise independent. It is also complete even if $\Pi$ is not.    

It is easily seen that we have.

\textit{For every pair $A, B\in \mathscr U_i$, if $A\neq B$, then $A$ and $B$ belong to the same partition $\pi\in \Pi_i$ if and only if $\hat{A}$, $\hat{B}$ are disjoint.}

Let $X^k=Y_1\times\cdots\times Y_k\times X_{k+1}\times\cdots\times X_n$. For $C\in \mathscr B$, let us define $C^k=\hat{C}_1\times\cdots\times \hat{C}_k\times C_{k+1}\times\cdots \times C_n$. If $\mathscr H$ is a subfamily of $\mathscr B$ , then $\mathscr H^k=\{C^k\colon  \text{for $C \in H$}\}$. 

Let $G$ be a polybox. One can wish to describe the corresponding mapping $G\mapsto G^k$ for polyboxes. Suppose $\mathscr G\subseteq \mathscr B$ is a suit for $G$. One can expect that $G^k$ can be defined as the union of $\mathscr G^k$. Then one has to check that this union is independent of the choice of the suit $\mathscr G$. This is done by induction. We prove only the first step, as the subsequent steps are mere repetitions. 

For every $x=(x_2, \ldots , x_d)\in X_2\times \cdots\times X_d$ we define two lines $\ell_x=X_{1}\times\{ x_2\}\times\cdots\times\{x_d\}$ and 
$\ell^1_x=Y_{1}\times \{x_2\}\times\cdots\times\{x_d\}$. Suppose first that $G\cap \ell_{x}=\ell_{x}$. Let $\mathscr G_x=\{ A\in \mathscr G\colon A\cap \ell_x\neq \emptyset\}$. It is a Keller family as a subfamily of $\mathscr G$. Therefore, there is $\pi \in\Pi_1$ such that $\pi =\{A_1\colon \text{for $A$ in $\mathscr G_x$}\}$. Then $\hat{\pi}$ is a partition of $Y_1$ and $\{\ell^1_x\cap \hat{A}\colon \text{for $A \in \mathscr G_x$}\}$ is a partition of $\ell^1_x$. Thus $G^1\cap \ell^1_x=\ell^1_x$ which  is evidently independent of the choice of $\mathscr G$. Suppose now that $G\cap \ell_x\neq\ell_x$.  Then, by Lemma \ref{C forte}, there is a unique $\pi \in \Pi_1$ such that $\rho=\{A_1\colon \text{for $A \in \mathscr G_x$}\}$ is a proper subfamily of $\pi$. In fact, $\pi$ is determined by the union of $\rho$ which coincides with the image of $G\cap \ell_x$ under the projection of $X$ onto $X_1$. Therefore, the set $\hat{\rho}=\{\hat {A}_1\colon \text{for $A_1\in \rho$}\}$ is independent of the choice of $\mathscr G$. And since the image of $G^1\cap\ell^1_x$ under the projection of $X^1$ onto $Y_1$ coincides with the union of $\hat{\rho}$, the set $G^1\cap\ell^1_x$ is also independent of the choice of $\mathscr G$. 

For $C\in \mathscr B$, $\mathscr G\subseteq \mathscr B$, and a polybox $G$ let us set: $\hat C= C^d$, $\hat {\mathscr G}={\mathscr G}^d$, $\hat G = G^d$. We can summarize all what we have shown as follows:

\begin{itemize}
\item[$(\gamma_1)$] \textit{For every pair $C, D\in \mathscr B$, if $C\neq D$, then $\{C, D\}$ satisfies Keller's condition if and only if $\hat{C}$, $\hat{D}$ are disjoint.}

\item[$(\gamma_2)$]  \textit{If $\mathscr G_1$ and $\mathscr G_2$ are suits for a polybox $G$, then $\hat {\mathscr G}_1$ and $\hat{\mathscr G}_2$ are suits for $\hat G$.} 
\end{itemize}

Now, we can complete the proof of Theorem \ref{D} easily.  Since  every $G\in \mathscr G$ is a proper box and for every $i\in \{1,\ldots, d\}$ and every nontrivial $\pi \in \Pi_i$ the cardinality of $\pi$ is $n$, it follows  by the definition of $\hat {G}$ that $|\hat {G}|= \frac 1{n^d} |Y|$.  As $\mathscr G$  and $\{X\}$ are Keller partitions of $X$, it follows that $\hat {\mathscr G}$  is a (Keller) partition of $\hat{X}=Y$.  Clearly, $|\mathscr G|= |\hat{ \mathscr G}|$.  Therefore,
$$
|Y|=\sum_{G\in \mathscr G} |\hat{G}|=  |\mathscr G| \frac {|Y|}{n^d} ,
$$
 which  implies (\ref{N}).
 $\square$

\section{Concluding remarks}

\textbf{1}. Theorem \ref{C} adresses a particular kind of tori: $\boldsymbol m=(n,\ldots, n)$. It is of some interest to extend this result to arbitray $\boldsymbol m=(m_1, \ldots, m_d)$. The following question arises:
\begin{quote}
Suppose $p=p(T)$ attains its maximum for a tiling $\mathscr G=[0,1)^d\oplus T$ of $\mathbb T^d_{\boldsymbol m}$. Is it true that $\mathscr G$ must be a multipile?
\end{quote}
One can easily check that if the tiling $\mathscr G$ is a multipile, then there is an ordering $i_1,i_2,\ldots, i_d$ of the set $\{1,\ldots, d\}$ such that 
$$
p(T)=1+ m_{i_1}+m_{i_1}m_{i_2}+\cdots +m_{i_1}m_{i_2}\cdots m_{i_{d-1}}.
$$
Therefore, if the answer to our question is affirmative and the ordering $i_1,i_2,\ldots, i_d$ is chosen so that $m_{i_1}\geq m_{i_2}\geq\ldots \geq m_{i_{d}}$, then the latter equality expresses the maximum value of $p$.  

\medskip
\noindent\textbf{2}.  To convince the reader  of the usefulness of the construction presented in the preceding section, we offer the following immediate consequence of $(\gamma_1)$ and $(\gamma_2)$:      
\begin{quotation}\textit{
Let $\mathscr G_1, \ldots, \mathscr G_s$  and $\mathscr H_1, \ldots, \mathscr H_s$ be two pairwise disjoint families that are contained in a system of boxes $\mathscr B$. Suppose that  for each $i$ the families  $\mathscr G_i$ and $\mathscr H_i$ are suits for the same polybox.  Moreover, suppose that $\mathscr G_1\cup\cdots \cup \mathscr G_s$ is  a suit for a polybox $G$.   Then $\mathscr H_1\cup\cdots\cup H_s$ is a suit for  $G$ as well.
} 
\end{quotation}  
Similar constructions have been widely exploited in \cite[Proposition 2.6, Section 10]{KP1}. 

\medskip
\noindent\textbf{3}.  Kisielewicz \cite{Ki} is concerned with the problem of finding unit cube partitions $[0,1)^d\oplus T$ of $\mathbb T_{\boldsymbol m}^d$, where $\boldsymbol m=(2,\ldots,2)$, for which the parameters $p_i(T)$ are equal and  as large as possible.  The simplest version of Theorem \ref{C} is already mentioned there \cite[Theorem 4.2, pp. 95, 104]{Ki}, and attributed to the present author.

\noindent
{Wydzia{\l} Matematyki, Informatyki i Ekonometrii, Uniwersytet Zielonog\'orski}\\
{ul. Podg\'orna 50, 65-246 Zielona G\'ora, Poland}\\
{\sf K.Przeslawski@wmie.uz.zgora.pl}

\end{document}